\documentstyle[12pt]{article}
\newcommand{\ltimes}{\bowtie} 
\newcommand{\bfTheta}{\bf{\Theta}} 
\newcommand{\Ut}{\hat{U}_t(g(A, \bfTheta))} 
\newcommand{\Uq}{U_q(g(A, \bfTheta))} 
\newcommand{\U}{U(g(A, \bfTheta))}

\newcommand{\ba}{\begin{eqnarray}}
\newcommand{\na}{\end{eqnarray}}
\newcommand{\ban}{\begin{eqnarray*}}
\newcommand{\nan}{\end{eqnarray*}}

\begin{document} 
\title{\normalsize{\bf SYMMETRIZABLE QUANTUM AFFINE SUPERALGEBRAS\\  
AND THEIR REPRESENTATIONS  }}  
\author{\small R. B. ZHANG\\  
\small  Department of Pure Mathematics\\  
\small  University of Adelaide\\ 
\small  Adelaide, S. A.,  Australia}    
\date{ \small January 1996  }
\maketitle
 
\begin{abstract} 
Aspects of the algebraic structure and representation theory of the 
quantum affine superalgebras with symmetrizable Cartan matrices 
are studied. The irreducible integrable highest weight representations 
are classified, and shown to be deformations of their classical 
counterparts. It is also shown that Jimbo type quantum affine 
superalgebras can be obtained by deforming universal 
enveloping algebras of ordinary (i.e., non-graded) affine algebras 
supplemented by certain parity operators.
{\bf Running title}: Quantum affine superalgebras  
\end{abstract} 

\pagebreak 
\normalsize

%\vspace{3cm}

\noindent
{\bf I.  Introduction} 
 
Quantum affine superalgebras are of great 
importance for the study of supersymmetric integrable models in 
statistical mechanics and quantum field theory. 
Recent research has also indicated that such algebraic structures 
may play a significant role in characterizing vacua  of $4$ - 
dimensional supersymmetric Yang - Mills theories and string 
compactifications.  Apart from their physical applications, 
quantum affine superalgebras are interesting from a  
mathematical point of view as well. They have many similarities  
to the ordinary (i.e., non - graded ) 
quantum affine algebras, thus a thorough investigation of  
their structures should be possible. 
It is also hoped that a representation theory can be developed for them, 
which will be workable in applications. 
One can quantize the affine superalgebras following Drinfeld and Jimbo 
relatively  easily,  once a proper understanding of the Serre type of 
presentations at the classical level is achieved. 
However, much more effort seems to  be  
required in order to develop their representation theory, as  
there already exist severe difficulties at the classical level.  
Although various special results are known in the area, e.g., 
the classification of the finite dimensional irreducible representations 
of $U_q(\widehat{gl}^{(1)}(m|n))$, there has been no attempt to study the 
quantum affine superalgebras systematically.  

The aim of this note is to investigate the structure and representation 
theory of the quantum affine superalgebras with symmetrizable Cartan 
matrices.  Their classical counterparts, which were classified by 
Kac\cite{Kac}, constitute the only class of affine superalgebras 
with  a well developed representation theory.  
One of our results is the classification of the irreducible integrable 
highest weight representations of these quantum affine superalgebras.     
We will generalize Lusztig's method \cite{Lusztig} to show that such 
representations are in one to one correspondence with 
the irreducible integrable highest weight representations  
of the associated classical affine superalgebras.   
Another result is that  quantum affine superalgebras can be obtained 
by deforming the universal enveloping algebras of {\em ordinary} 
affine algebras supplemented by certain parity operators, wherein  
some kind of Bose - Fermi transmutation is exhibited. 
This result will be useful physically, e.g., for showing equivalences 
of various integrable models. Mathematically, it also bears considerable 
implications in the classification of quantum affine superalgebras 
and representation theory.  In this note, we will use the result to show  
a correspondence between the representations of the super and ordinary  
quantum affine algebras.

The arrangement of the paper is as follows.  In section II we define 
the Drinfeld type of quantum affine superalgebras and examine 
some of their algebraic features from the point of view of deformation 
theory.  In section III  we classify the integrable highest weight irreps, 
and in section IV we study the afore mentioned Bose - Fermi transmutation.\\

\noindent 
{\bf II.  Quantum Affine Superalgebras}

Let $A=\left(A_{i j}\right)_{i, j=0}^n$ be the Cartan matrix of an 
affine Lie superalgebra which satisfies the following conditions
\ban a_{i i}=2,  & a_{i j}\le 0, \  \ i\ne j, \\ 
     a_{i j}=0 \  \  \mbox{iff}\   \  a_{j i}=0, 
     & a_{i j}\in 2{\bf Z}, \ \ \mbox{if} \ \  i\in \bfTheta,   
\nan 
where $\bfTheta$ is a nonempty subset of the index set $I=\{0, 1, ..., n\}$. 
Such Cartan matrices are called symmetrizable, 
and the affine Lie superalgebras associated with them  
have been classified. They are given by the following 
Dynkin diagrams.     
   
\vspace{.5cm}  
\begin{center}  Table available upon request  \end{center}  
\vspace{.5cm}  

\noindent 
In the above table,  a diagram has $n+1$ nodes with the $i$ - th node 
being white if $i\not\in\bfTheta$, and black if $i\in\bfTheta$.  The $i$ - th 
and $j$ - th nodes are connected by $\mbox{max}(|a_{i j}|,  |a_{j i}|)$ 
lines; if $|a_{i j}|>|a_{j i}|$, the lines are endowed with an 
arrow pointing towards the $i$ - th node.  The numerical marks for the 
diagram will be denoted by $a_i$,  $i=0, 1, ..., n$, which satisfy the 
condition $\sum_{j=0}^n a_{i j} a_j =0$.

We denote by $g(A, \bfTheta)$  the complex affine superalgebra 
associated with the Cartan matrix $A$ and the subset $\bfTheta\subset I$. 
Let $H^*$ be the dual vector space of the  Cartan subalgebra of 
$g(A, \bfTheta)$. 
Then $H^*$ has a basis $\{\Lambda_0,  \alpha_i, \ i\in I\}$, where the 
$\alpha_i$ are the simple roots.  A 
nondegenerate bilinear form $(. \ , . )$ on $H^*$ can be defined  
in the standard way,  satisfying,   
\ban 
{{2(\alpha_i, \alpha_j)}\over{(\alpha_i, \alpha_i)} } = a_{ i j}, 
\ \ \   
{{2(\Lambda_0, \alpha_i)}\over{(\alpha_i, \alpha_i)} } = \delta_{ i 0},
\ \ \ 
(\Lambda_0, \Lambda_0)=0. 
\nan 
An appropriate normalization for the form can always be chosen 
such that  
\ban 
(\alpha_\mu, \ \alpha_\mu)&=&1,   \ \ \ \forall \mu\in\bfTheta, 
\nan 
and we will work with this normalization throughout.

The Drinfeld model of quantum affine superalgebra $\Ut$ 
is a ${\bf Z}_2$ graded associative algebra over the 
ring ${\bf C}[[t]]$, with $q=exp ( t )$, completed 
with respect to the $t$ - adic topology of ${\bf C}[[t]]$. 
It is generated by the elements 
$\{ d, \ h_i, \ e_i, \ f_i, \ i\in I\}$,  
subject to the relations 
\ba 
%\begin{array}{ll} 
k_i k^{-1}_i=1, &   k_i k_j = k_j k_i,  \nonumber \\
{[}d, k^{\pm 1}_i\} =0, & \nonumber \\ 
{[}d, e_i\}=\delta_{i 0} e_i, 
& [d, f_i\}=-\delta_{i 0} f_i,  \nonumber \\ 
k_i e_j = q^{(\alpha_i, \ \alpha_j)} e_j k_i, 
&k_i f_j= q^{-(\alpha_i, \ \alpha_j)} f_j k_i,\nonumber \\ 
{[}e_i, f_j\} = \delta_{i j}{{k_i - k_i^{-1}}\over
{q^{\epsilon_i} - q^{-\epsilon_i} }}, 
&\forall i, j\in I, \nonumber \\ 
(Ad e_i)^{1-a_{i j}}(e_j) = 0, &(Ad f_i)^{1-a_{i j}}(f_j)=0, 
\ \   \forall i\ne j,  
%\end{array} 
\label{definition} 
\na    
where  
\ban 
k_i &=& q^{h_i}, \\ 
\epsilon_i&=&\left\{\begin{array}{l l} 
                  1, & \mbox{if}\  (\alpha_i, \alpha_i)=1,\\ 
                  1, & \mbox{if}\  (\alpha_i, \alpha_i)=2,\\ 
                  2, & \mbox{if}\  (\alpha_i, \alpha_i)=4.  
               \end{array}\right.  
\nan   
All the generators are chosen to be homogeneous, 
with $d$, $h_i, \ i\in I$, and $e_j$, $f_j$, $j\not\in \bfTheta$,  
being even, and $e_\mu$, $f_\mu$, $\mu\in\bfTheta$, being odd.  
For a homogeneous element $x$, we define $[x]=0$ if $x$ is even, 
and $[x]=1$ when odd.  The graded commutator $[ .\  ,\ . \}$  
represents the usual commutator when any one of the two arguments 
is even, and the anti commutator when both arguments are odd.
The adjoint operation $Ad$ is defined by
\ban 
Ad e_i(x)&=& e_i x - (-1)^{[e_i][x]} k_i x k_i^{-1} e_i, \\ 
Ad f_i(x)&=& f_i x - (-1)^{[f_i][x]} k_i^{-1} x k_i f_i. 
\nan      
For $x$ being a monomial in $e_j$'s or $f_j$'s, it carries a 
definite weight $\omega(x)\in H^*$. Then 
$Ad e_i(x)$ $=$ $e_i x $ $-$ 
$(-1)^{[e_i][x]} q^{(\alpha_i, \ \omega(x))} x e_i$, 
and similarly for $Ad f_i(x)$.

The quantum affine superalgebra $\Ut$ has the structures of a 
${\bf Z}_2$ graded Hopf algebra with a co - multiplication 
\ban 
\Delta(d)&=&d\otimes 1 + 1\otimes d, \\     
\Delta(h_i)&=&h_i\otimes 1 + 1\otimes h_i, \\     
\Delta(e_i)&=&e_i\otimes k_i + 1\otimes e_i, \\     
\Delta(f_i)&=&f_i\otimes 1 + k_i^{-1}\otimes h_i.    
\nan 
A co - unit and an antipode also exit, but we shall not spell 
them out explicitly, as they will not be used here. 
Our main concern in this letter is the algebraic structures 
and the representations of the quantum affine superalgebras.

An important fact is that  as a ${\bf Z}_2$ graded associative algebra,  
$\Ut$ is a deformation of the universal enveloping algebra 
$\U$ of $g(A, \bfTheta)$ in the sense of \cite{Gerstenhaber},  that is, 
being a topologically free ${\bf C}[[t]]$ module, $\Ut$ is isomorphic 
to the ${\bf C}[[t]]$ module $\U[[t]]$, consisting of power series 
in $t$ with coefficients in $\U$,  and there also exists the  
algebra isomorphism $\Ut/t\Ut \cong \U$.  
This of course is a standard fact in the theory of quantum 
groups \cite{Drinfeld} \cite{Flato}.  However, proving it is a highly 
nontrivial matter, and is well out of the scope of this letter. 

To make things more explicit, 
let $m$ be the associative multiplication of $\U$.     
Denote by $m_t$ the associative multiplication of $\Ut$. 
Then $m_t$  is a ${\bf C}[[t]]$  bi - linear map 
$m_t: \U[[t]]{\hat\otimes}\U[[t]]$ $\rightarrow \U[[t]]$ 
of the form $m_t = m$  $+$ $ \sum_{i=1}^\infty t^i m^{(i)}$, 
where $m$ is the multiplication of $\U$, and 
$\hat\otimes$ is the tensor product completed with respect to 
the $t$ - adic topology of ${\bf C}[[t]]$.  The $m^{(i)}:$ 
$\U$ $\otimes$ $\U\rightarrow \U$ are ${\bf Z}_2$ graded vector 
space maps, which are homogeneous of degree zero.  
Associativity of $m_t$ imposes stringent conditions on the
maps $m^{(i)}$. In particular, the first nonvanishing $m^{(i)}$ 
must be a $2$ - cocycle in the language of Hochschild cohomology.  
In view of the fact that  the Drinfeld quantum affine algebras 
are nontrivial deformations of the universal enveloping algebras 
of the associated affine algebras, we expect the deformations defining 
the quantum affine superalgebras also to be nontrivial.

Consider a $\Ut$ module $V_t$, with the module action 
denoted by  $\circ_t$. If $V_t$ is a free ${\bf C}[[t]]$ module, 
then $V_t=V[[t]]$, with $V=V_t/t V_t$ a complex vector space. 
Assume that for any given $a\in\U\subset\Ut$, $v\in V\subset V_t$, 
\ban 
a\circ_t v &=& a\circ v + o(t)\in V[[t]],
\nan 
where $\circ$ represents a ${\bf C}$ bi - linear map 
$\U\otimes V \rightarrow V$, then $\circ$ defines a module action 
of $\U$ on $V$. To see that our claim is indeed correct,
consider another element $ b \in $ $\U$ $\subset\Ut$. Then
\ban
b\circ_t[ a\circ_t (v + t V_t)]&=& m_t(b, a)\circ_t v + t V_t\\
     &=& m(b, a) \circ v + t V_t.
\nan

Conversely,  let the complex vector space $V$ be a $\U$ module, 
with the module action $\circ$.  
If there exists a ${\bf C}[[t]]$ bi - linear map
$\circ_t:  \U[[t]]{\hat \otimes} V[[t]] \rightarrow V[[t]]$, such that
for any $a, b\in\U\subset\Ut$,  $v\in V\subset V[[t]]$, 
\ban 
a\circ_t v &=& a\circ v + o(t)\in V[[t]],\\
a\circ_t(b\circ_t v)&=&m_t(a, b)\circ_t v,
\nan
then $V[[t]]$ furnishes a $\Ut$ module. In this case, we say that 
the $\Ut$ module $(V[[t]], \circ_t)$ is a deformation of the $\U$ 
module $(V, \circ)$, and the representation of $\Ut$ afforded by 
$(V[[t]], \circ_t)$ the deformation of the representation 
of $\U$ furnished by $(V, \circ)$. (Note the difference between our 
definition of deformation of modules and that of \cite{Pinczon}. )   
We will call the deformation {\em trivial} if there exists a
${\bf C}[[t]]$ linear map $\Phi_t = id + t\phi_1 + t^2\phi_2 + ... :$ 
$\U[[t]]\rightarrow \U[[t]]$ such that $ a\circ_t v = \Phi_t(a)\circ v$,
$\forall a\in \Ut$, $v\in V[[t]]$, where $\circ$ is 
${\bf C}[[t]]$ - linearly extended to $\Ut$. 
Needless to say, not all representations of $\U$ can 
be deformed into representations of $\Ut$. It is a very interesting 
problem to characterize the deformability of a $\U$ module in 
cohomological terms, and we hope to return to the problem in 
the future.  We should also mention that if an irreducible 
representation can be deformed at all, then the deformation must be 
trivial.\\ 

\noindent 
{\bf III. Integrable Highest Weight Modules} 
 
Let us first construct the irreducible highest weight $\Ut$ modules.   
We will omit the symbols $m_t$ and $\circ_t$ from our notations 
whenever confusion is not likely to arise. 
Let $U^+_q$ be the ${\bf Z}_2$ graded subalgebra of $\Ut$ generated by 
the $h_i$, $e_i$, $i=0, 1, ..., n$, together with $d$, 
and $N^-_q$ that generated by the $f_i$, $i=0, 1, ..., n$. 
Let $v^\Lambda_+\otimes {\bf C}[[t]]$ be a one dimensional 
$U^+_q$ module satisfying 
\ban 
h_i v^\Lambda_+ =(\Lambda, \alpha_i)\  v^\Lambda_+,  
& e_i v^\Lambda_+ =0,  & \forall i=0, 1, ..., n,\\ 
&d v^\Lambda_+ = 0. &  
\nan 
We construct the ${\bf C}[[t]]$ module 
\ban 
{\overline V}_t(\Lambda)&=& \Ut\otimes_{U^+_q} v^\Lambda_+, 
\nan 
which is clearly isomorphic to $N^-_q\otimes v^\Lambda_+$, and therefore, 
is spanned by  the elements of the form  
$f_{i_1} f_{i_2} ... f_{i_p}\otimes v^\Lambda_+$, $i_s\in I$, $p\in{\bf Z}_+$. 
(We will ommit the tensor product sign $\otimes$ from such expressions 
hereafter.)   
Define a bi - linear  action of $\Ut$ on  ${\overline V}_t(\Lambda)$ by  
\ba 
d\ (f_{i_1} f_{i_2} ... f_{i_p} v^\Lambda_+)
&=& -\sum_{s=1}^p \delta_{i_s 0}\  f_{i_1} f_{i_2} ... f_{i_p} v^\Lambda_+, 
\nonumber \\ 
k_i \ (f_{i_1} f_{i_2} ... f_{i_p} v^\Lambda_+) 
&=& q^{(\Lambda - \sum_{s=1}^p \alpha_{i_s},\   \alpha_i)} 
\ f_{i_1} f_{i_2} ... f_{i_p} v^\Lambda_+, \nonumber\\ 
f_i\ (f_{i_1} f_{i_2} ... f_{i_p} v^\Lambda_+) 
&=& f_i f_{i_1} f_{i_2} ... f_{i_p} v^\Lambda_+,\nonumber \\   
e_i\ (f_{i_1} f_{i_2} ... f_{i_p} v^\Lambda_+)
&=&\sum_{s=1}^p \delta_{i i_s} (-1)^{[e_{i}]\sum_{k=1}^{s-1} [f_{i_k}]} 
\ f_{i_1} f_{i_2} ... {\hat f}_{i_s}... f_{i_p} v^\Lambda_+ \nonumber  \\    
&\times& 
{ {q^{(\Lambda - \sum_{r=s+1}^p \alpha_{i_r},\   \alpha_{i_s})} - 
q^{- (\Lambda - \sum_{r=s+1}^p \alpha_{i_r},\  \alpha_{i_s})} } 
\over{q^{\epsilon_{i_s}} - q^{-\epsilon_{i_s}} }}.  \label{module}  
\na  
All the relations of (\ref{definition}) are clearly satisfied, 
except the Serre relations amongst the $e_i$'s. 
Set $S_{i j}=(Ad e_i)^{1-a_{i j}}(e_j)$, $i\ne j$, 
It is a consequence of the `quadratic' 
relations that $[S_{i j}, \ f_k\}=0$. Thus for all $i\ne j$,  
\ban 
S_{i j} (f_{i_1} f_{i_2} ... f_{i_p} v^\Lambda_+) &=&0,   
\nan  
and ${\overline V}_t(\Lambda)$ indeed yields a $\Ut$ module. 

This module is in general not irreducible, but contains a maximal 
proper submodule $M(\Lambda)$ such that 
\ban 
V_t(\Lambda)&=&{\overline V}_t(\Lambda)/M(\Lambda), 
\nan  
yields an irreducible $\Uq$ module, which is called 
an irreducible highest weight module with highest weight $\Lambda$. 
The image of $v^\Lambda_+$ under the canonical projection is 
the maximal vector of $V_t(\Lambda)$. Standard arguments show that 
up to isomorphisms, $V_t(\Lambda)$ is uniquely determined by 
its highest weight.

Following the terminology of the representation theory of 
Lie algebras, we call a $\Ut$ module $V_t$ integrable if all 
$e_i$ and $f_i$ act on $V_t$ by locally nilpotent endomorphisms, 
namely, for any $v\in V_t$, there exists a nonnegative integer 
$m_v<\infty$ such that 
\ban 
(e_i)^{m_v} v = &(f_i)^{m_v} v =0, & \forall i\in I. 
\nan

Consider the  irreducible highest weight $\Ut$ module $V_t(\Lambda)$ 
with highest weight $\Lambda$ and maximal vector $v^\Lambda_+$. 
It is obviously true that the $e_i$ 
always act on $V_t(\Lambda)$ by locally nilpotent endomorphisms.  
However, nilpotency of the $f_i$ action imposes strong  
conditions on the highest weight.  

For a fixed $i\not\in\bfTheta$, the elements $e_i$, $f_i$ and $h_i$ 
generate a $U_{q^{\epsilon_i}}(sl(2))$ subalgebra of $\Ut$.  In order for 
$(f_i)^{m} v^\Lambda_+$ to vanish, $\Lambda$ must satisfy the condition 
${2{(\Lambda, \ \alpha_i)}\over{(\alpha_i, \ \alpha_i)}}\in {\bf Z}_+$. 
When $\mu\in\bfTheta$, we have normalized 
$(\alpha_\mu, \ \alpha_\mu)=1$.  Now $e=e_\mu$, $f=f_\mu$, $h=h_\mu$, 
generate a $U_q(osp(1|2))$ subalgebra, 
\ban 
[h, \ e]= e, & [h, \ f]= -f, & e f + f e = { { q^h - q^{-h}} 
                       \over{q - q^{-1}}}.   
\nan        
Since $f^{m} v^\Lambda_+=0$ for a large enough $m$, 
but $v^\Lambda_+\ne 0$, there must exist an integer $k$, 
$0< k < m$, such that $f^k v^\Lambda_+ \ne 0$, 
and  $f^{k+1} v^\Lambda_+=0$. 
Applying $e$ to $f^{k+1} v^\Lambda_+$, we arrive at 
\ban 
e f^{k+1} v^\Lambda_+ &={ { q^{{k+1}\over{2}} - q^{-{{k+1}\over{2}}}}
\over{(q - q^{-1})(q^{1/2} - q^{-1/2})} } 
\left[ q^{{{(\Lambda, \ \alpha_\mu)}\over{(\alpha_\mu, \ \alpha_\mu)}} 
- {{k}\over{2}}}
- (-1)^k q^{-{{(\Lambda, \ \alpha_\mu)}\over{(\alpha_\mu, \ \alpha_\mu)}}
+{{k}\over{2}}} \right]  f^k v^\Lambda_+ &=0,     
\nan  
which requires ${{2(\Lambda, \ \alpha_\mu)}\over{(\alpha_\mu, \ \alpha_\mu)}}
=k$, and $k\in 2{\bf Z}_+$. In fact,\\ 
{\em The irreducible highest weight $\Ut$ module $V_t(\Lambda)$ with 
highest weight $\Lambda$ is integrable if and only if 
\ba
{2{(\Lambda, \ \alpha_i)}\over{(\alpha_i, \ \alpha_i)}}\in {\bf Z}_+, 
&\forall i\in I, \nonumber\\  
{{2(\Lambda, \ \alpha_\mu)}\over{(\alpha_\mu, \ \alpha_\mu)}}\in 2{\bf Z}_+, 
&\forall \mu\in\bfTheta. \label{integral}  
\na } 

Note the presence of  the second condition requiring 
the Dynkin labels associated with the odd simple roots be 
non - negative {\em even} integers, which is not needed in the case of 
ordinary quantum affine algebras.   
We prove the assertion  following the strategy of
\cite{Lusztig}.  As pointed out earlier, all the $e_i$ act on 
$V_t(\Lambda)$ by locally nilpotent endomorphisms. We have also 
seen that under the given conditions of $\Lambda$, the maximal 
vector $v^\Lambda_+$ of $V_t(\Lambda)$ is annihilated by  
a sufficiently high power of each $f_i$, $i\in I$. 
Now consider the element $w=f_{i_1} f_{i_2} ... f_{i_p} v^\Lambda_+$. 
We use induction on $p$ to prove the nilpotency of the 
action of the $f_i$ on $w$. Assume that $x=f_{i_2} ... f_{i_p} v^\Lambda_+$ 
is annihilated by $(f_i)^m$, $\forall i\in I$. 
Then $(f_{i_1})^m w = (f_{i_1})^{m+1} x=0$. 
For $j\ne i_1$, consider $(f_j)^{m-a_{j i_1}} w = 
(f_j)^{m-a_{j i_1}} f_{i_1} x$.  
By using the Serre relation  $(Ad f_j)^{1-a_{j i_1}}$  
$f_{i_1}$ $=0$, we can express $(f_j)^{m-a_{j i_1}} f_{i_1}$ as 
a ${\bf C}[[t]]$ linear combination of the elements 
$(f_j)^{-a_{j i_1}- \nu} f_{i_1} (f_j)^{m+\nu}$,  
$\nu=0, 1, ..., -a_{j i_1}$,  which all annihilate $x$.     
Hence, $(f_j)^{m-a_{j i_1}} w =0$.

As $V_t(\Lambda)$ can be generated by repeatedly applying 
the $f_i$ to the maximal vector $v^\Lambda_+$,  a subset $\cal B$ 
of all the elements of the form $f_{i_1} f_{i_2} ... f_{i_p} v^\Lambda_+$ 
provides a basis of $V_t(\Lambda)$ over ${\bf C}[[t]]$. 
Hence we have proved that all $e_i$ and $f_i$ act on $V_t(\Lambda)$ 
by locally nilpotent endomorphisms, and $V_t(\Lambda)$ is integrable.

Let $V(\Lambda)$ be the vector space over ${\bf C}$ with the basis 
$\cal B$.  Then as a ${\bf C}[[t]]$ module, $V_t(\Lambda)=V(\Lambda)[[t]]$. 
Our earlier discussions assert that $V(\Lambda)\cong $ 
$V_t(\Lambda)/t V_t(\Lambda)$ carries a natural $\U$ module
structure, and $V_t(\Lambda)$ is a deformation of $V(\Lambda)$.       
It follows the integrability of $V_t(\Lambda)$ that 
$V(\Lambda)$ is integrable as a $\U$ module. It is also of 
highest weight type, and is cyclically generated by a the maximal vector 
with weight $\Lambda$.  A result of \cite{Kac} states that 
an integrable $\U$ module is completely reducible. 
Thus $V(\Lambda)$, being cyclically generated, must be irreducible. 
Also recall that every integrable irreducible highest weight $\U$ module 
is uniquely determined by an element $\Lambda\in H^*$ satisfying 
the same conditions as (\ref{integral}). Thus, \\  
{\em Every irreducible integrable highest weight $\Ut$ module is a  
deformation of an irreducible integrable highest weight $\U$ module  
with the same highest weight,  and all such irreducible $\U$ modules 
can be deformed. }

Note that the integrable lowest weight irreps of the quantum affine 
superalgebras can be studied in the same way, and the above result 
applies as well.   It should also be mentioned that Kac' 
character formula \cite{Kac} for the integrable
highest weight irreps of $\U$ still works in the quantum case.\\ 

\noindent
{\bf IV.  Jimbo Model and Bose -  Fermi Transmutation} 
 
The Jimbo version of quantum affine superalgebra $U_q(g(A, \bfTheta))$
is a ${\bf Z}_2$ graded associative algebra over the complex
number field ${\bf C}$, generated by the elements
$\{ d,\ k_i, \ k_i^{-1},$ $e_i$, $f_i, \ i\in I\}$,
subject to the same relations as (\ref{definition}), 
but with $q$ now being regarded as a non - zero complex parameter. 
Nevertheless,  it is possible to formulate the Jimbo type 
`quantization' within the framework of deformation theory 
\cite{Flato}.

We still set $q=exp(t)$, and let $t_i =t\epsilon_i$. 
Define
\ban
S_i&=& { { k_i - k_i^{-1}}\over{q^{\epsilon_i} - q^{-\epsilon_i}} },\\ 
C_i&=& { { k_i + k_i^{-1}}\over{2} }.  
\nan
The relations of (\ref{definition}) involving $k_i^{\pm 1}$ 
can now be re - expressed in terms of $S_i$ and $C_i$, 
\ban 
[d, C_i]=0, & & [d, S_i]=0, \\  
C_i S_j = S_j C_i, &  &   
(C_i)^2  - (S_i)^2\  sinh^2 t_i = 1, \\  
C_i e_j - e_j C_i \ cosh[t(\alpha_i, \ \alpha_j)] 
       &=& e_j S_i \ sinh t_i \ sinh[t(\alpha_i, \ \alpha_j)], \\ 
S_i e_j - e_j S_i \ cosh[t(\alpha_i, \ \alpha_j)] 
        &=& e_j C_i \ sinh[t(\alpha_i, \ \alpha_j)]/sinh t_i, \\  
C_i f_j - f_j C_i \ cosh[t(\alpha_i, \ \alpha_j)]  
       &=&- f_j S_i \ sinh t_i \ sinh[t(\alpha_i, \ \alpha_j)], \\ 
S_i f_j - f_j S_i \ cosh[t(\alpha_i, \ \alpha_j)]  
        &=& - f_j C_i \ sinh[t(\alpha_i, \ \alpha_j)]/sinh t_i, \\ 
{[}e_i, \ f_j\}&=& \delta_{i j} S_i,   
\nan   
while the Serre relations remain the same. 
We can now regard $\Uq$ as generated by $d, \ C_i, \ S_i$,  $e_i, \ f_i$  
for any $t\in {\bf C}$. Furthermore, for a fixed $t_0\in {\bf C}$, 
and $t=t_0 +\tau $, we can consider $\tau$ as a formal parameter, 
and define the formal Jimbo quantum affine superalgebra $\Uq$ as 
a properly completed ${\bf C}[[\tau]]$ algebra generated by $d, \ C_i$, 
$S_i, \ e_i, \ f_i$ with the same relations. Then $\Uq$ is  a
deformation of $U_{exp(t_0)}(g(A, \bfTheta))$. 

At $t_0=0$, $U_{exp(t_0)}(g(A, \bfTheta))$ is isomorphic to an 
extension of the universal enveloping algebra of $g(A, \bfTheta)$ 
by the $C_i$, which satisfy $C_i^2=1$. Explicitly, 
\ban 
U_1(g(A, \bfTheta))&=&\U\otimes {\bf C Z}_2^{\otimes (n+1)},  
\nan 
where ${\bf C Z}_2^{\otimes (n+1)}$ is the group algebra of the
abelian group generated by the $C_i$. Therefore, strictly speaking, 
the Jimbo model of  $\Uq$ is  not a deformation of $\U$, 
but rather an an extension of $\U$ by some parity operators.

More interesting is the case when $t_0=i\pi$. 
At $\tau=0$, the relations become 
\ban
[d, C_i]=0, & & [d, S_i]=0, \\  
C_i S_j = S_j C_i, &  & (C_i)^2   = 1, \\
C_i e_j =  (-1)^{(\alpha_i, \ \alpha_j)} \ e_j C_i, 
& &C_i f_j = (-1)^{(\alpha_i, \ \alpha_j)}\ f_j C_i, \\
S_i e_j - (-1)^{(\alpha_i, \ \alpha_j)} \ e_j S_i 
&=&(-1)^{(\alpha_i, \ \alpha_j)+\epsilon_i} 
\ (\alpha_i, \ \alpha_j)\ e_j C_i, \\
S_i f_j -(-1)^{(\alpha_i, \ \alpha_j)}\ f_j S_i 
&=&- (-1)^{(\alpha_i, \ \alpha_j)+\epsilon_i} 
\ (\alpha_i, \ \alpha_j)\ f_j C_i, \\
{[}e_i, \ f_j\}&=& \delta_{i j} S_i.  
\nan
and the Serre relations read
\ban 
(ad e_i)^{1-a_{i j}}(e_j) = 0, &(ad f_i)^{1-a_{i j}}(f_j)=0, 
&i\ne j,  
\nan 
with  
\ban 
ad e_i (x)&=& e_i x - (-1)^{(\alpha_i, \ \omega(x))+[x][e_i]} x e_i, \\ 
ad f_i (x)&=& f_i x - (-1)^{(\alpha_i, \ \omega(x))+[x][e_i]} x f_i. 
\nan 

These are the defining relations for the complex associative algebra 
$U_{-1}(g(A, \bfTheta))$, which, unfortunately, are rather complicated, 
and not very illuminating. However, by applying certain inner automorphisms 
constructed out of the $C_i$, we can cast  the relations into a more 
familiar form. 

For definiteness, let us consider $B^{(1)}(0, n)$.    
Set $\sigma_i=\prod_{k=i}^n C_k$, and define 
\ban
D&=&d, \\  
H_i&=&(-1)^{\epsilon_i} C_i\ S_i,\\ 
E_i&=&\sigma_{i+1} \ (\sigma_1)^{\delta_{i 0}}\ e_i, \\  
F_i&=&\sigma_i\ (\sigma_1)^{\delta_{i 0}} \ f_i.  
\nan 
Now something rather intriguing happens: these elements do {\em not} 
obey the defining relations of the affine superalgebra $B^{(1)}(0, n)$, 
instead they generate the universal enveloping algebra of the twisted 
ordinary ( i.e., non - graded)  affine Lie algebra $A^{(2)}_{2 n}$. 
Recall that the universal enveloping algebra of $A^{(2)}_{2 n}$ and 
that of $B^{(1)}(0, n)$ are totally different algebraic structures, 
although their underlying vector spaces  ( ignoring the ${\bf Z}_2$ 
grading in the case of $B^{(1)}(0, n)$ ) are isomorphic.  
Nevertheless, $U_q(B^{(1)}(0, n))$ can be obtained as a 
deformation of the universal enveloping algebra of 
$A^{(2)}_{2 n}$ supplemented by $n+1$ parity operators: 
a kind of transmutation between the ordinary 
affine algebra (which is boson - like)
and affine superalgebra ( which is fermion - like)
takes places upon quantization. 
Such a transmutation was found in the case of $osp(1|2 n)$ and 
$so(2 n +1)$ in \cite{Zhang}, providing a natural explanation  
for the observation made by Rittenberg and Scheunert \cite{Scheunert} 
that there was a one to one correspondence between the tensorial 
irreducible representations of $so(2 n +1)$ and the finite dimensional 
irreducible representations of $osp(1|2 n)$.  

Note that the $C_i$ generate the group algebra of the abelian group 
${\bf Z}_2^{\otimes (n+1)}$. When acting by conjugation on the 
elements of the $A^{(2)}_{2 n}$ generators, they give rise to 
parity factors,  i.e., $\pm$  signs.   We introduce the notation 
$U(A^{(2)}_{2 n})\ltimes {\bf C Z}_2^{\otimes (n+1)}$ to 
illustrate the fact that $U_{-1}(B^{(1)}(0, n))$ is 
the universal enveloping algebra of $A^{(2)}_{2 n}$ 
supplemented by the $C_i$.

A case by case study shows that such Bose - Fermi transmutation occurs
with other affine superalgebras as  well;  
we have  
\ba 
U_{-1}(B^{(1)}(0, n))&\cong& U(A^{(2)}_{2 n})\ltimes 
          {\bf C Z}_2^{\otimes (n+1)}, \ \ n>1,\nonumber\\  
U_{-1}(B^{(1)}(0, 1))&\cong& U(A^{(2)}_{2 })\ltimes 
          {\bf C Z}_2^{\otimes 2}, \nonumber\\
U_{-1}(A^{(2)}(0, 2n-1))&\cong& U(B^{(1)}_n)\ltimes 
           {\bf C Z}_2^{\otimes (n+1)}, \ \ n>2, \nonumber\\
U_{-1}(A^{(2)}(0, 3))&\cong& U(C^{(1)}_2)\ltimes
           {\bf C Z}_2^{\otimes 3}, \nonumber\\
U_{-1}(C^{(2)}(n+1))&\cong& U(D^{(2)}_{n+1})\ltimes
           {\bf C Z}_2^{\otimes (n+1)}, \nonumber\\
U_{-1}(C^{(2)}(2))&\cong& U(A^{(1)}_1)\ltimes
          {\bf C Z}_2^{\otimes 2}. \label{transmutation}  
\na 
However, the $A^{(4)}(0, 2n)$ series proves to be  an exception   
\ba
U_{-1}(A^{(4)}(0, 2n))&\cong&U(A^{(4)}(0, 2n))\ltimes
           {\bf C Z}_2^{\otimes (n+1)}, \ \ n=1, 2, ...,  
\na 
where no Bose - Fermi transmutation has been observed.

The transmutation between ordinary quantum affine algebras
and quantum affine superalgebras can also be realized at 
the level of representations. 
Consider an irreducible integrable highest weight 
module $V_t(\Lambda)$ of the Drinfeld quantum affine superalgebra 
$\Ut$ studied in the last section.  Note that $t$ enters 
the formulae (\ref{module}) through $q$, 
thus by specializing $t$ to a complex number 
$t=i\pi + \tau$, with $\tau$ a generic complex parameter, 
we obtain from $V_t(\Lambda)$ a module of the Jimbo quantum affine 
superalgebra $\Uq$, which we denote by ${\check V}_q(\Lambda)$.  
If $g(A, \bfTheta)$ is one of the affine superalgebras appearing 
in (\ref{transmutation}), then the representation of $\Uq$ furnished 
by ${\check V}_q(\Lambda)$ can be 
realized by an irreducible integrable highest weight 
representation of the ordinary quantum affine algebra 
$U_{-q}(g(A, \emptyset))$, where $g(A, \emptyset)$, appearing on the 
right hand sides of (\ref{transmutation}),  is the ordinary affine Lie 
algebra with the same Cartan matrix $A$, but with all generators 
being even.

For the sake of concreteness, consider again the case 
of $U_q(B^{(1)}(0, 2 n))$.  Denote by $D, \ E_i, \ F_i, \ (K_i)^{\pm 1}$ the 
generators of $U_{-q}(A^{(2)}_{2 n})$, while the generators of 
$U_q(B^{(1)}(0, 2 n))$ are still denoted by $d, \ e_i, \ f_i, \ (k_i)^{\pm 1}$. 
Since the Cartan subalgebras of $B^{(1)}(0, 2 n)$ and $A^{(2)}_{2 n}$ are
isomorphic, we will make no distinctions between them. 

Let ${\check W}_{-q}(\Lambda)$ be an irreducible $U_{-q}(A^{(2)}_{2 n})$ 
module with highest weight $\Lambda$ satisfying the conditions 
(\ref{integral}).  As a complex vector space ${\check W}_{-q}(\Lambda)$ 
admits the weight space decomposition
\ban 
{\check W}_{-q}(\Lambda)&=&\bigoplus_{\omega\le\Lambda} W^{\omega}, 
\nan 
where each $W^{\omega}$ is finite dimensional, and $(\alpha_i, \ \omega)
\in {\bf Z}$, $\forall i\in I$. 
Define a $U_q(B^{(1)}(0, 2 n))$ action on ${\check W}_{-q}(\Lambda)$ by 
\ban 
d  w &=& D w, \\  
k_i w&=& (-1)^{(\alpha_i, \ \omega)} K_i w,\\  
e_i w &=& (-1)^{(\beta_{i+1}-\beta_1 \delta_{i 0}, 
            \ \omega + \alpha_i)} E_i w, \\ 
f_i w&=& (-1)^{(\beta_i - \beta_1 \delta_{i 0}, 
            \ \omega + \alpha_i)} F_i w,    
            \ \ \ \ \ \ \forall w\in W^{\omega},   
\nan
where $\beta_{i}=\sum_{r=i}^n\alpha_r$.  
Direct calculations show that this definition 
indeed preserves the defining relations of $U_q(B^{(1)}(0, 2 n))$, 
thus turning ${\check W}_{-q}(\Lambda)$ into a $U_q(B^{(1)}(0, 2 n))$ 
module. This module is clearly irreducible, and has highest weight 
$\Lambda$. Thus it is isomorphic to ${\check V}_{q}(\Lambda)$.  
Observe that the subset of $H^*$ satisfying (\ref{integral}) exhausts all 
the integral dominant weights for $B^{(1)}(0, 2 n)$. 
Therefore, every irreducible integrable highest weight 
representation of $U_q(B^{(1)}(0, 2 n))$  can be realized this way. 

In a similar way we can show that the same result also holds for  
other affine superalgebras:  \\ 
{\em Let $g(A, \bfTheta)$ be an affine superalgebra appearing in 
(\ref{transmutation}). Then each irreducible integrable highest weight 
representation of $\Uq$ can be realized by a representation of 
$U_{-q}(g(A,\ \emptyset))$ of the same kind.}

However, the converse is not true. There exist integrable irreps of 
$U_{-q}(g(A,\ \emptyset))$ with highest weights not satisfying the 
second condition of (\ref{integral}).  
 
%\pagebreak  


\vspace{2cm} 
  
      
\begin{thebibliography}{9999} 
\bibitem{Kac}     V. G. Kac, Adv. Math. {\bf 30} (1978) 85. 
\bibitem{Lusztig} G. Lusztig, Adv. Math. {\bf 70} (1988) 237. 
\bibitem{Gerstenhaber}M. Gerstenhaber, Ann. Math. {\bf 78} (1963) 267.
\bibitem{Drinfeld} V. G. Drinfeld,  Leningrad Math. J.  {\bf 1} (1989) 321. 
\bibitem{Flato} P. Bonneau, M. Flato, M. Gerstenhaber and 
                  G. Pinczon, Commun. Math. Phys., {\bf 161} (1994) 125. 
\bibitem{Pinczon} M. Lesimple and  G. Pinczon, J. Math. Phys., 
                  {\bf 34} (1993) 4251.  
\bibitem{Zhang}   R. B. Zhang,  Lett. Math. Phys., {\bf 25} (1992) 317.    
\bibitem{Scheunert} V. Rittenberg and M. Scheunert, 
                    Commun. Math. Phys., {\bf 83} (1982) 1. 
\end{thebibliography}
\end{document}